\documentclass[reqno,11pt]{amsart}

\usepackage{graphicx}
\usepackage{url}
\usepackage[cmtip,arrow]{xy}  
\usepackage{pb-diagram,pb-xy}  
\usepackage{graphicx}    
\usepackage{amsfonts}
\usepackage{amsmath, amssymb}
\usepackage[latin1]{inputenc}
\usepackage{color, enumerate}

\usepackage{pgfplots}
\usepackage{subfigure}
\usepackage{caption}
\usepackage{float}

\usepackage{amsmath}
\usepackage{amsthm}
\usepackage{amsfonts}
\usepackage{array}
\newcolumntype{C}[1]{>{\centering\arraybackslash}p{#1}}

\newcommand{\T}{\mathbb{T}}

\newcommand{\UU}{\mathbb{U}}

\newcommand{\C}{\mathbb{C}}

\newcommand{\R}{\mathbb{R}}

\newcommand{\without}[1]{\backslash\{#1\}}

\newcommand{\pp}[2]{\frac{\partial#1}{\partial#2}}

\newcommand{\norm}[1]{\|#1\|}

\newtheorem{proposition}{Proposition}
\newtheorem{theorem}[proposition]{Theorem}
\newtheorem{definition}[proposition]{Definition}

\theoremstyle{remark}
\newtheorem{remark}[proposition]{Remark}

\begin{document}

\bibliographystyle{alpha}

\title[Hamiltonian systems in singular symplectic manifolds]{Examples of integrable and non-integrable systems on singular symplectic manifolds}

\author{Amadeu Delshams}
\address{Departament de Matem\`{a}tiques,
 Universitat Polit\`{e}cnica de Catalunya, ETSEIB, Avinguda Diagonal 647, Barcelona}\email{amadeu.delshams@upc.edu}
\author{Anna Kiesenhofer}
\address{Departament de Matem\`{a}tiques, Universitat Polit\`ecnica de Catalunya, EPSEB, Avinguda del Doctor Mara\~{n}\'{o}n 44--50, Barcelona}
\email{anna.kiesenhofer@upc.edu}
\author{Eva Miranda}
\address{Departament de Matem\`{a}tiques,
 Universitat Polit\`{e}cnica de Catalunya and Barcelona Graduate School of Mathematics BGSMath, EPSEB, Avinguda del Doctor Mara\~{n}\'{o}n 44--50, Barcelona}
\email{eva.miranda@upc.edu}
\thanks{Amadeu Delshams is partially supported by the Catalan grant 2014SGR504. Anna Kiesenhofer is supported by AGAUR FI doctoral grant. Anna Kiesenhofer and Eva Miranda are partially supported by MINECO grant MTM2012-38122-C03-01/FEDER and by the Catalan Grant 2014SGR634.}

\date{\today}
\begin{abstract}
We present a collection of examples borrowed from celestial mechanics and projective dynamics. In these examples symplectic structures with singularities arise naturally from regularization transformations, Appell's transformation  or classical changes  like McGehee coordinates, which end up blowing up the symplectic structure or lowering its rank at certain points. The resulting geometrical structures that model these examples are no longer symplectic but symplectic with singularities which are mainly of two types: $b^m$-symplectic and $m$-folded symplectic structures. These examples comprise the three body problem as non-integrable exponent and some integrable reincarnations such as the two fixed-center problem. Given that the geometrical and dynamical properties of $b^m$-symplectic manifolds  and folded symplectic manifolds are well-understood \cite{guimipi, guimipi2,gmps, kms, martinet, folded, Gualtierili, gualtierietal, marcutosorno2, scott, gmw}, we envisage that this new point of view in this collection of examples can shed some light on classical long-standing problems concerning the study of dynamical properties of these systems seen from the Poisson viewpoint.
\end{abstract}
\maketitle


\section{Introduction}
Integrability and non-integrability of some classical problems in physics and celestial mechanics such as the Kepler problem of the $2$ or $3$-body problems is well-understood \cite{arnoldkozlov}. Even if the $3$-body problem is not integrable  some restricted cases like the 2-fixed center problem are integrable. When studying such systems ad hoc transformations have been considered in order to understand their dynamics (for instance in the McGehee change of coordinates) and integrability (Appell's transformation for Newton's systems).

In this article we provide a list of classical examples in celestial mechanics and projective dynamics and analyze the classical changes done in the theory to study their geometrical and dynamical properties. We observe that these classical changes induce singularities in the Darboux symplectic structure of the phase space (cotangent bundle). Transformations that preserve the Darboux symplectic form are classically known as canonical. Some of these singularities have a physical interpretation and correspond to actual collisions of the bodies in the $n$-body problem or correspond to lines at infinity.

Classically there have been two different approaches in the literature concerning these examples: One approach, in the vein of  \cite{mcgehee, moeckel, moeckel2, devaney},  is built on implementing the transformation, realizing that these transformations are not canonical, forgoing the geometrical structures and just focusing on the computational simplifications that these changes  yield. A different approach represented by a group including  \cite{celletti, bolotin, gmm}  is based on the idea of implementing the transformations in a canonical way. This typically implies cumbersome computations so as to change the momenta accordingly. Some of these papers (cf. \cite{gmm}) adhere to this point of view because standard KAM theory is employed. Some progress understanding the KAM theory in the singular setting has already been achieved in \cite{kms}. This new point of view could be helpful in this scenario.

In this paper we propose  an intermediate approach which explores and imports the best of both worlds.

The singularities induced by these classical transformations are geometrically well-understood from the point of view of symplectic and Poisson geometry. They are either $b^m$-symplectic structures, $m$-folded symplectic structures or some variants of these. By identifying the geometrical structure behind them, we can benefit from the understanding of their dynamics attained in \cite{guimipi, guimipi2, gmps, kms, martinet, folded}. In particular, new results concerning the existence of periodic trajectories which has been a center of attention in the theory (see for instance \cite{chenciner}) could be unraveled with this new perspective in mind.

\section{Singular symplectic structures}\label{sec:pre}

As it is typical in physics, a Poisson structure is regarded as a bracket on functions in the following way:
$$\{f,g\}:= \Pi (df, dg),\qquad f,g\in C^\infty(M)$$ \noindent where $\Pi$ is a bivector field.

In this section we define $b$-Poisson/$b$-symplectic and folded symplectic structures as well as $b^m$-symplectic structures and $m$-folded symplectic structures. The rest of structures present in this article will be either a sophistication of these or a combination of both types (like the ones which show up in projective dynamics which can be understood in the Dirac framework).

The set of singular structures which we consider in this paper share a common feature: they are symplectic away from a singular set (hypersurface) and they behave \lq\lq nicely" on this surface. These structures have a simple
Darboux canonical form  at points on the critical set which we state below.

A Poisson structure of $b$-type $\Pi$ on a manifold $M^{2n}$ defines a symplectic structure on a dense set in $M^{2n}$. The set of points where the Poisson structure is not symplectic is a hypersurface  defined as the vanishing set of the multivectorfield $\Pi^n$. Namely,

\begin{definition}
Let $(M^{2n},\Pi)$ be an oriented Poisson manifold such that the map
\begin{equation}\label{bpoisson}
p\in M\mapsto(\Pi(p))^n\in\Lambda^{2n}(TM)
\end{equation}
is transverse to the zero section, then $Z=\{p\in M|(\Pi(p))^n=0\}$ is a hypersurface and we say that $\Pi$ is a \textbf{$b$-Poisson structure} on $(M^{2n},Z)$ and $(M^{2n},Z)$ is a \textbf{$b$-Poisson manifold}. The hypersurface $Z$ is called \textbf{critical hypersurface}.
\end{definition}

\begin{remark}
Observe that the transversality condition in the definition above is equivalent to $0$ being a regular value of the map $p\in M\mapsto(\Pi(p))^n\in\Lambda^{2n}(TM)$.
\end{remark}

These Poisson structures were classified in dimension $2$ by Radko \cite{radko}. Their higher dimensional study was motivated by the works of Melrose \cite{melrose} on the calculus on manifolds with boundary and by the works of Nest and Tsygan \cite{nestandtsygan} on deformation quantization of manifolds with boundary. Its systematic study started with the work of Guillemin, Miranda and Pires \cite{guimipi, guimipi2} and has attracted the attention of other mathematicians interested in the geometry and topology of these manifolds (see for instance \cite{gualtierietal, frejlichmartinezmiranda, marcutosorno2, gmps, kms, planas}).

The geometrical and dynamical properties of $b$-Poisson manifolds is well-understood thanks to the fact that it is possible to address problems on these manifolds using generalized De Rham forms which are called $b$-forms which admit \emph{poles} at the critical hypersurface $Z$. Grosso modo, these forms are just dual to the multivectorfields which have maximal rank away from the critical hypersurface. This is an important tool in the theory which allows to import many successful techniques in the symplectic realm such as Moser's path method. In this paper we will work with this dual approach.

In particular a Darboux theorem holds for these manifolds,
\begin{theorem}[\textbf{$b$-Darboux theorem, \cite{guimipi2}}]\label{theorem:bDarboux}
Let $\omega$ be a $b$-symplectic form on $(M^{2n},Z)$ and $p\in Z$. Then we can find a local coordinate chart $(x_1,y_1,\ldots,x_n,y_n)$ centered at $p$ such that locally the hypersurface $Z$ is locally defined by $y_1=0$ and
$$\omega=d x_1\wedge\frac{d y_1}{y_1}+\sum_{i=2}^n d x_i\wedge d y_i.$$
\end{theorem}

The dual $b$-Darboux theorem gives a local normal form of type,

\begin{equation}
\Pi = y_1\frac{\partial}{\partial x_1}\wedge
\frac{\partial}{\partial y_1}+\sum_{i=2}^n \frac{\partial}{\partial x_i}\wedge
\frac{\partial}{\partial y_i}.
\end{equation}

It is possible to generalize these structure and consider more degenerate singularities of the Poisson structure. This is the case of $b^m$-Poisson structures \cite{scott} for which $\omega^m$ has a singularity of $A_m$-type in Arnold's list of simple singularities \cite{arnold0, arnold}. In the same spirit we may consider other singularities in this list.

As it happens with $b$-Poisson structures, it is possible and convenient  to consider a dual approach in their study and work with forms. We refer the reader to \cite{scott} and \cite{gmw} for details.
Recall,

\begin{definition} A {\bf symplectic $b^m$-manifold}  is a pair $(M^{2n}, Z)$ with a closed $ b^m $-two form $ \omega $  which has maximal
 rank  at every $p \in M$.
\end{definition}

In \cite{gmw} a $b^m$-Darboux theorem is proved for $b^m$-Poisson structures,

\begin{theorem}[\textbf{$b^m$-Darboux theorem, \cite{gmw}}]\label{theorem:bnDarboux}
Let $\omega$ be a $b^m$-symplectic form on $(M^{2n},Z)$ and $p\in Z$. Then we can find a coordinate chart $(x_1,y_1,\ldots,x_n,y_n)$ centered at
$p$ such that the hypersurface $Z$ is locally defined by $y_1=0$ and
$$\omega=d x_1\wedge\frac{d y_1}{y_1^m}+\sum_{i=2}^n d x_i\wedge d y_i.$$
\end{theorem}

In the same way, dually we obtain a $b^m$-Darboux form,

\begin{equation}
\Pi = y_1^m\frac{\partial}{\partial x_1}\wedge
\frac{\partial}{\partial y_1}+\sum_{i=2}^n \frac{\partial}{\partial x_i}\wedge
\frac{\partial}{\partial y_i}
\end{equation}

A second class of important geometrical structures that model some problems in celestial mechanics  are \emph{folded symplectic structures}. These are closed $2$-forms on even dimensional manifolds which are non-degenerate on a dense set thanks to the following transversality condition.

\begin{definition}
Let $(M^{2n},\omega)$ be a manifold with $\omega$ a closed $2$-form such that the map
$$p\in M\mapsto(\omega(p))^n\in\Lambda^{2n}(T^*M)$$
is transverse to the zero section, then $Z=\{p\in M|(\omega(p))^n=0\}$ is a hypersurface and we say that $\omega$ defines a \textbf{folded symplectic structure } on $(M,Z)$ and $(M,Z)$ is a \textbf{folded symplectic manifold}. The hypersurface $Z$ is called \textbf{folding hypersurface}.
\end{definition}

The normal forms of folded symplectic structures was studied by Martinet \cite{martinet}.

 \begin{theorem}[\textbf{folded-Darboux theorem, \cite{martinet}}]\label{theorem:foldedDarboux}
Let $\omega$ be a folded symplectic form on $(M^{2n},Z)$ and $p\in Z$. Then we can find a local coordinate chart $(x_1,y_1,\ldots,x_n,y_n)$ centered at
$p$ such that  the hypersurface $Z$ is locally defined by $y_1=0$ and
$$\omega=y_1d x_1\wedge{d y_1}+\sum_{i=2}^n d x_i\wedge d y_i.$$
\end{theorem}

 In analogy to the case of $b^m$-symplectic structures we define a new class of folded structures, namely \textbf{$m$-folded symplectic structures} for which $\omega^n$ has singularities of $A_m$-type in Arnold's list of simple singularities \cite{arnold0}.

 For them $\omega^n$ has a local normal form of type  $\omega^n= y_1^m dx_1\wedge\dots\wedge dy_n$

Integrability in this context can be easily understood as integrability in the symplectic dense set. That is to say in the examples considered in this paper we require our $2n$-dimensional manifold which is a cotangent bundle with a form with poles or zeros to admit a set of $n$ functionally independent commuting functions with respect to the induced bracket. This notion coincides with the one for Poisson structures defined in \cite{camilleevapol} and has also been considered in \cite{kms}.
\section{The Kepler problem and Levi-Civita transformation}
 In this section  we consider several classical changes  for the $n$-body problem, typically non-symplectic. First these transformations are considered to solve the $2$-body problem (which is integrable) and second to study the systems close to singularities like collisions or line at infinity.
\subsection{The two-body problem}

The two-body problem is the system consisting of two bodies with masses $m_1, m_2$ and positions $q_1,q_2 \in R^3, q_1 \neq q_2$ moving under their mutual gravitational attraction. Hence according to Newton's law of gravity the equations of motion are
\begin{align}\label{newton}
m_i \ddot{q_i} = \mathcal{G} m_1 m_2 \frac{ q_j-q_i}{\norm{q_2-q_1}^3}= \pp{U}{q_i}, \quad i,j=1,2, \quad i\neq j,
\end{align}
where $\mathcal{G}$ is the gravitational constant and we have introduced the negative gravitational potential
$$ U:= m_1 m_2 \frac{\mathcal{G}}{\norm{q_2-q_1}}.$$

To obtain a Hamiltonian formulation of the problem we define the momenta $p_i = m_i \dot q_i$ and the Hamiltonian
$$H:= E_{kin} - U = \frac{\norm{p_1}^2}{2m_1}+\frac{\norm{p_2}^2}{2m_2} - U.$$
Then the second-order differential equations \eqref{newton} are equivalent to the first-order systems $(i=1,2)$
\begin{align*}
\dot q_i = \pp{H}{p_i}, \qquad \dot p_i = -\pp{H}{q_i},
\end{align*}
which is Hamiltonian with respect to the canonical symplectic form
$$\omega_0 = d q_1 \wedge d p_1 + d q_2 \wedge d p_2 .$$

\subsection{The Kepler problem - Jacobi coordinates}

By introducing an appropriate set of coordinates (see for instance \cite{meyerhalloffin}) we can transform the two-body problem to the problem of one body moving in a central force field (Kepler problem). Let
\begin{align}
g &= \nu_1 q_1 + \nu_2 q_2, \qquad &G &= p_1 + p_2,\\
w &= q_2-q_1,				\qquad &W & = -\nu_2 p_1 + \nu_1 p_2,
\end{align}
where $\nu_i = m_i/(m_1+m_2)$. Note that $g$ is the center of mass and $G$ is the total linear momentum. The coordinate $w$ is the relative position of the second body with respect to the first one. The other ``momentum" coordinate $W$ is chosen in such a way that the change of coordinates is {\it canonical} (i.e., the symplectic form is preserved). The coordinates  $(g,w,G,W)$ are called Jacobi coordinates.

In these coordinates the Hamiltonian is
$$H(g,w,G,W)=\frac{\norm{G}^2}{2 \nu} + \frac{\norm{W}^2}{2 M}-\mathcal{G} \frac{m_1 m_2}{\norm{w}}$$
where $\nu = m_1 + m_2$ and $M=m_1 m_2/(m_1 + m_2)$.

Writing down the Hamiltonian equations explicitly
\begin{align*}
\dot g &= \pp{H}{G} = \frac{G}{\nu}, \qquad &\dot G &=-\pp{H}{g}=0,\\
\dot w &= \pp{H}{W} = \frac{W}{M}, \qquad &\dot W &=-\pp{H}{w} = -\frac{m_1 m_2 w}{\norm{w}^3},
\end{align*}
we see that total linear momentum $G$ is preserved and that the center of mass moves with constant velocity  $\frac{G}{\nu}$. Hence the problem reduces to the second line of equations; physically this means that we are viewing the system from the perspective of one body with coordinates $w$ under the influence of the central force field of a body with mass $M$. The upshot is that we are dealing with a Hamiltonian system on $(\R^3\backslash\{0\})\times \R^3$ with Hamiltonian function
$$H(w,W)=\frac{\norm{W}^2}{2 M}-\mathcal{G} \frac{m_1 m_2}{\norm{w}}.$$
This is precisely what is known as the {\bf Kepler problem}.

\subsection{The planar Kepler problem and Levi-Civita coordinates}

The Kepler problem has three degrees of freedom. An easy calculation shows that the angular momentum $w\times W$
is a conserved quantity. Hence the position and momentum vector lie in the same plane; we can consider the system as {\it planar}. After performing a rotation we can assume that the relevant components of $w$ and $W$ are the first two. We discard the third component and are left with a system on $\R^2\backslash\{0\} \times \R^2$, the {\bf planar Kepler problem} with Hamiltonian
$$H(w,W)=\frac{\norm{W}^2}{2 M}-\mathcal{G} \frac{m_1 m_2}{\norm{w}}, \qquad (w,W)\in (\R^2 \backslash \{0\}) \times \R^2.$$
The corresponding second-order differential equation is
\begin{equation}\label{planarkepler}
\ddot w = -\mathcal{G} \frac{(m_1 + m_2)w}{\norm{w}^3}, \qquad w \in \R^2.
\end{equation}

\subsection{Levi-Civita coordinates}

We identify $\R^2$ with $\C$ and view Equation \eqref{planarkepler} as a differential equation over $\C$. In this way we can introduce a new coordinate $u$ via
$$\frac{u^2}{2} = w.$$
This step is called the {\bf Levi-Civita regularization procedure}.

\begin{remark}The physical motivation for doing this transformation is that the Kepler problem takes on a very simple form:  Performing a change of time $dt = r d\tau$ and fixing an energy orbit $H=h$,  the second-order differential equation \eqref{planarkepler} becomes
$$u'' + \omega u = 0 ,$$
the well-known harmonic oscillator. Here, the prime denotes differentiation with respect to $\tau$ and $\omega = \sqrt{h/ 2}$.
\end{remark}


\subsection{Geometrical structure}

We want to study how the symplectic form changes under the Levi-Civita regularization. (The change of time is not of interest for us.)

Recall that the phase space coordinates are $(w,W)\in  (\R^2 \backslash \{0\}) \times \R^2$ and the symplectic form is the canonical one. We leave the momentum coordinates $W$ unchanged and express the space coordinates $w$ in terms of $u$:
$$w_1+iw_1 = w = \frac{u^2}{2} = \frac{(u_1 + i u_2)^2}{2} = \frac{u_1^2}{2}-\frac{u_2^2}{2} +iu_1u_2.$$
The real resp. imaginary part correspond to the components $w_1,w_2$. Hence their differentials are
$$dw_1 = u_1 du_1 -  u_2 du_2, \qquad dw_2 = u_1 du_2 + u_2 du_1$$
and the symplectic form is
\begin{align*}
\omega = dw_1 \wedge dW_1 + dw_2 \wedge dW_2 = &\  u_1 du_1 \wedge dW_1 - u_2 du_2 \wedge dW_1 \\
&+ u_1 du_2 \wedge dW_2 + u_2 du_1\wedge dW_2.
\end{align*}

Observe that $\omega\wedge\omega=(u_1^2-u_2^2)   du_1 \wedge dW_1\wedge du_2 \wedge dW_2$. So $\omega\wedge\omega$ does not cut the zero section of $\Lambda^2 T^*M$ transversally (which is the condition for a form to be folded \cite{folded}) so $0$ is not a regular value of the function $\omega\wedge\omega$  but  the singularity is  non-degenerate of hyperbolic type.

\begin{remark}
A different version of the Levi-Civita transformation described in \cite{bolotin, celletti} does not keep the momenta unchanged but transforms them in such a way that the total change is {\it symplectic}. More precisely the change is given by:
$$ w = \frac{u^2}{2}, \qquad W= \frac{U}{\overline u},$$
where $(u,U)$ are the new coordinates. In this case the symplectic form remains the standard one but the equations look more involved.
\end{remark}

\subsection{The Kepler problem in three dimensions - KS transformation}

We now skip the step of reducing the Kepler problem to a planar system and consider Equation \eqref{planarkepler} in three dimensions. This is relevant e.g. for studying binary collisions in the three-body problem (which cannot be restricted to two dimensions in general).

 Instead of working with complex numbers, the regularization procedure now employs the {\bf quaternion algebra} $\UU$, see \cite{waldvogel}. Recall that $\UU$ consists of objects of the form
$$u=u_0 + i u_1 + j u_2 + k u_3$$
where $i,j,k$ are the three independent "imaginary" units and multiplication is defined via the non-commutative laws
$$ij = -ji = k,  jk = -kj = i, ki = -ik = j.$$
We identify the quaternion $u$ with the vector $u=(u_0,u_1,u_2,u_3)\in \R^4$. Moreover, we introduce the star conjugation
$$u^\ast := u_0 + i u_1 + j u_2 - k u_3$$
Now instead of considering the squaring of complex numbers as in the previous section, we define the mapping
\begin{equation}\label{kstransform}
u \mapsto w:= \frac{ u u^\ast }{2}.
\end{equation}
The image is the set of quaternions with vanishing $k$ component and can be identified with $\R^3$. The preimage of a number $w=v v^\ast$ is given by the one-parameter family of quaternions of the form $v \cdot e^{k \theta} := v \cdot (\cos \theta + k \sin \theta )$ where $\theta\in S^1$.


Writing the transformation \eqref{kstransform} explicitly (where we identify the image with $\R^3$), we have
\begin{align*}
 w_0 &= \frac{ u_0^2 - u_1^2 - u_2^2 + u_3^2 }{2}\\
 w_1 & = u_0 u_1-u_2 u_3 \\
 w_2 & = u_0 u_2 + u_1 u_3.
\end{align*}
This is known as the {\bf Kustaanheimo-Stiefel (KS)} transformation. We choose a solution with vanishing $k$-component, i.e. $u_3 = 0$. Denoting the original space coordinates by $(w_0,w_1,w_2)$ with conjugate momenta $(W_0,W_1,W_2)$, the symplectic form becomes in the new space coordinates $(u_0,u_1,u_2)\in \R^3$:
\begin{align*}
\omega =& dw_0\wedge dW_0 + dw_1 \wedge dW_1 + dw_2 \wedge dW_2 = \\
 =& (u_0du_0-u_1 du_1- u_2 du_2) \wedge dW_0 + (u_0du_1+u_1du_0) \wedge dW_1 \\&+ (u_0 du_2 + u_2 du_0)\wedge dW_2.
 \end{align*}

In order to determine which kind of geometric structure this defines, we compute  $\omega\wedge\omega\wedge\omega= (u_0^3-u_1^2u_0-u_2^2 u_0) du_0\wedge dW_0\wedge du_1\wedge dW_1 \wedge du_2\wedge dW_2$. Observe that the coefficient of $\omega\wedge\omega\wedge\omega$ is $2u_0 w_0$. This is a sophistication of $m$-folded symplectic structures.

\section{Total collapse  in the $n$-body problem}
In this section we focus on triple collisions for the three-body problem but the computation holds \emph{mutatis mutandis} for the total collapse in the $n$-body problem.

Following \cite{moeckel, moeckel2} and \cite{mcgehee}, consider the system of three bodies with masses $m_1,m_2, m_3$ and positions $\mathbf{q_1}=(q_1,q_2,q_3),\mathbf{q_2}=(q_4,q_5,q_6), \mathbf{q_3}=(q_7,q_8,q_9)\in \R^3$. Similarly we denote the components of the momenta by $p_1,\ldots,p_9$. We define the $9\times 9$ matrix $M:=\textup{diag}(m_1,m_1,m_1,m_2,m_2,m_2,m_3,m_3,m_3)$.

We assume that we are working in central coordinates, so the centre of mass remains at the origin:
$$m_1 \mathbf{q_1} + m_2 \mathbf{q_2}+m_3 \mathbf{q_3}=0.$$

We introduce the following ``McGehee"-coordinates:

\begin{equation}\label{triple}
r:= \sqrt{q^T M q},\qquad s:=\frac{q}{r}, \qquad z:=p\sqrt{r}.
\end{equation}

Note that $r=0$ corresponds to triple collisions. Essentially, these are spherical coordinates since $s$ lies on the unit-sphere in $R^9$ with respect to the metric given by $M$.

To specify a well-defined chart based on the formulas above \eqref{triple}, we restrict to the subset $q_9>0$ of the phase space $\R^9\times\R^9$ and consider the coordinates
$$(r,s_1,\ldots,s_8,z_1,\ldots,z_9).$$
The inverse of this chart is
\begin{align*}
q_i &= r s_1, \qquad i=1,\ldots ,8 \\
q_9 &= r\sqrt{ \frac{1-\sum_{i=1}^8 s_i^2 m_i}{m_9} }, \\
p_i &= \frac{z_i}{\sqrt{r}}, \qquad i=1,\ldots ,9. \\
\end{align*}
Computing the differentials one sees that the standard symplectic form $\sum_{i=1}^9 dq_i\wedge dp_i$ becomes
\begin{align*}
 \sum_{i=1}^8 &\left( \frac{s_i}{\sqrt{r}} dr \wedge dz_i + \sqrt{r} ds_i \wedge dz_i - \frac{z_i}{2 \sqrt{r}} ds_i \wedge dr \right) + \\
 +&\frac{1}{\sqrt{m_9 r \mu}} \left( \mu dr \wedge dz_9 - r \sum_{i=1}^8 m_i s_i ds_i \wedge dz_9 + \frac{z_9}{2} \sum_{i=1}^8 m_i s_i ds_i \wedge dr \right),
 \end{align*}
 where we have introduced the function $\mu:=1-\sum_{i=1}^8 s_i^2 m_i$ to simplify notation.

We compute the top wedge of the structure:
$$\bigwedge_{i=1}^9 dq_i \wedge dp_i = \sqrt{\frac{\mu r^7}{ m_9}} ds_1 \wedge dz_1 \wedge ds_2 \wedge dz_2 \wedge \ldots \wedge ds_8 \wedge dz_8 \wedge dr \wedge dz_9, $$
hence for $r=0$ this expression vanishes to order $\frac{7}{2}$ and is a $\frac{7}{2}$-folded symplectic structure.
In the $n$-body problem we should get $m$-folded symplectic structure for a certain $m$.

\section{The elliptic restricted three-body problem: McGehee coordinates}

Let us now consider a special case of the three-body problem where one of the bodies is assumed to have negligible mass. E.g. this happens if the system is given by Sun, Jupiter and an asteroid. Then a useful approximation is to assume that the motion of the two heavy bodies, called primaries (here Sun and Jupiter), is independent of the small body, hence given by Kepler's law for the two-body problem. Moreover, we assume that all the three bodies move in a plane.

We are interested in the resulting dynamical system for the small body (the asteroid), which moves under the influence of the time-dependent gravitational potential of the primaries
$$U(q,t)= \frac{1-\mu}{|q-q_1|} + \frac{\mu}{|q-q_2|},$$
where we assume that the masses of the primaries are normalized and given by $\mu$ (resp. $1-\mu$); their time-dependent positions are  $q_1=q_1(t)$ (resp. $q_2=q_2(t)$). \\

The Hamiltonian of the system is given by
$$H(q,p,t)= p^2/2 - U(q,t),\quad (q,p) \in \R^2 \times \R^2,$$
where $p=\dot q$ is the momentum of the planet. \\

The primaries are assumed to move around their center of mass on ellipses ({\it elliptic} restricted three-body problem). As explained in \cite{delshams}, it is useful to introduce polar coordinates to describe the motion of the small body; then for $q=(X,Y)\in \R^2\without{0}$, we have
$$ X = r \cos \alpha, Y = r \sin \alpha,\qquad  (r,\alpha) \in \R^+ \times \T$$
The momenta $p=(P_X, P_Y)$ are transformed in such a way that the total change of coordinates
$$(X,Y, P_X, P_Y) \mapsto (r,\alpha, P_r=:y, P_\alpha=:G)$$
is canonical, i.e. the symplectic structure remains the same.

To study the behaviour at $r=\infty$, it is standard to introduce the following {\bf McGehee coordinates}
$(x,\alpha,y,G)$, where $$ r=\frac{2}{x^2}, \quad x \in \R^{+}.$$

This transformation is non-canonical i.e. the symplectic structure changes and an easy calculation shows that for $x>0$ it is given by
$$-\frac{4}{x^3} dx \wedge dy + d\alpha \wedge d G.$$
This extends naturally to a $b^3$-symplectic structure  on $\R\times \T \times \R^2$ in the sense of \cite{scott}.

Equivalently, the Poisson bracket is
$$\{f,g\} = -\frac{x^3}{4} \left(\pp{f}{g}\pp{g}{y} - \pp{f}{y}\pp{g}{x}\right) + \pp{f}{\alpha}\pp{g}{G}-\pp{f}{G}\pp{g}{\alpha}$$
and such a Poisson structure is used extensively in \cite{delshams} to describe the dynamics close to the infinity manifold $x=y=0$.

\section{Examples from Projective Dynamics}

In this last section we examine an integrable reincarnation of a constrained $3$-body problem. The case of two fixed center-problem. This is a Newton system.

A Newton system is a system of the form: $\ddot{q}=f(q)$. Appell's transformation is a central projection which changes the \lq\lq screen" of
projection (hyperplane at infinity). Two such systems are equivalent and sometimes it is convenient to use these changes of screen in order to study some of their properties (such as integrability). This is the main principle of projective dynamics (see  \cite{albouy, albouy2}) which works well for Newton's systems.

 An outstanding example of these systems is the
 two fixed-center problem
(Euler, 1760). A particle in the plane moves under
the gravitational attraction of two fixed points $A$ and $B$ with
masses $m_A$ and $m_B$.

 This system can be written as:
{$$\ddot{q}=-m_A\frac{q_A}{\|{q_A}\|^3}-m_B\frac{q_B}{\|q_B\|^3}$$}
($q_A=q-A$, $q_B=q-B$).
Two first integrals are given by the Hamiltonian
$$H= \frac{1}{2}\|\dot{q}\|^2
-\frac{m_A}{\|{q_A}\|}-m_B\frac{m_B}{\|q_B\|}$$
and
$$G=\langle q_A\wedge \dot{q},q_B\wedge \dot{q}\rangle
-\frac{m_A}{\|{q_A}\|}\langle q_A,u\rangle-m_B\frac{m_B}{\|q_B\|}
\langle q_B,u\rangle$$
\noindent where $u=q_A-q_B$

The first integrals $H$ and $G$ Poisson commute thus defining an integrable system. In \cite{mirandaintegrable} following Albouy \cite{albouy} it is explained how to perform central projection for the two-center problem: Starting with the cotangent bundle in $T^*(\mathbb R^2)$
\lq\lq position" homogeneous coordinates are denoted by $[q_0 : q_1 :q_2]$. The initial affine chart is $q_0=1$ and a change to the affine chart $q_2=1$ is considered. Performing central projection to the  screen $q_2=1$, and changing the momenta  accordingly, this yields  an integrable system on  $q_2=1$.

The initial Darboux symplectic structure in $T^*(\mathbb R^2)$  becomes:
\begin{align*}
dv_1\wedge dq_1&+\frac{q_1}{q_2} (dq_1\wedge
dv_2+dq_2\wedge dv_1)+\\
 &+ \frac{(v_2q_1-v_1q_2)}{q_2^2}dq_1\wedge
dq_2+\left( \frac{q_1^2}{q_2^2}-1 \right) dv_2\wedge dq_2.
\end{align*}

This form has poles and zeros so for some hypersurfaces the structure is $b^m$-symplectic and for other hypersurfaces is $m$-folded symplectic. The right geometrical framework for these structures to coexist is that of Dirac structures  \cite{dirac}.

\end{document}